\input xy
\xyoption{all}

\documentclass{jhrs}
\usepackage{amsmath,latexsym,amssymb,mathrsfs}

\newtheorem{theorem}{Theorem}[section]

\newtheorem{corollary}[theorem]{Corollary}

\newtheorem{lemma}[theorem]{Lemma}

\newtheorem{remark}[theorem]{Remark}

\newcommand{\C}{\mathcal{C}}
\newcommand{\spaces}{\;\;\;}
\newcommand{\columnIIs}[2]{\fontsize{8}{9.2}\selectfont\left(\begin{array}{c}#1\\#2\end{array}\right)}
\newcommand{\matrixIIxIIs}[4]{\fontsize{8}{9.2}\selectfont\left(\begin{array}{cc}#1 & #2\\#3 & #4\\\end{array}\right)}
\newcommand{\ra}{\longrightarrow}
\newcommand{\T}{\mathcal{T}}
\newcommand{\x}{\times}
\newcommand{\V}{\mathcal{V}}

\def\theenumi{\alph{enumi}}

\begin{document}
\title{Closedness properties of internal relations IV: Expressing additivity of a category \\ via subtractivity}

\author{Zurab Janelidze}
\thanks{Partially supported by South African National Research Foundation.}
\email{zurab@gol.ge}
\address{Department of Mathematics and Applied Mathematics,
University of Cape Town, \\ Rondebosch 7701, Cape Town, South
Africa\\ \\ A. Razmadze Mathematical Institute of Georgian Academy
of Sciences,\\ 1 M.~Alexidze Street, 0193 Tbilisi, Georgia}


\bigskip
\begin{abstract}
The notion of a subtractive category, recently introduced by the
author, is a ``categorical version'' of the notion of a (pointed)
subtractive variety of universal algebras, due to A.\,Ursini. We
show that a subtractive variety $\C$, whose theory contains a unique
constant, is abelian (i.e.\! $\C$ is the variety of modules over a
fixed ring), if and only if the dual category $\C^\mathrm{op}$ of
 $\C$, is subtractive. More generally, we show that
$\C$ is additive if and only if both $\C$ and $\C^\mathrm{op}$ are
subtractive, where $\C$ is an arbitrary finitely complete pointed
category, with binary sums, and such that each morphism $f$ in $\C$
can be presented as a composite $f=me$, where $m$ is a monomorphism
and $e$ is an epimorphism.
\end{abstract}

\received{}   % receive date (for example: 11 October 1999)
\revised{}    % receive date
\published{}  % publish date
\submitted{Francis Borceux}  % Name of Journal's Editor, who submitted Article

\volumeyear{2006} % Volume Year
\volumenumber{1}  % Volume Number
\issuenumber{1}   % Issue Number

\startpage{1}     % PageNumber of first page

\maketitle
\section*{Introduction}

A variety $\V$ of universal algebras is \textit{subtractive} in the
sense of A.\,Ursini \cite{Urs94}, if the theory of $\V$ contains a
binary term $s$ (called a \textit{subtraction term}) and a nullary
term $0$, satisfying the identities $s(x,0)=x$ and $s(x,x)=0$. An
example of a subtractive variety is the variety of (additive)
groups, for which $s(x,y)=x-y$. More generally, any Mal'tsev variety
\cite{Mal54}, whose theory contains a nullary term, is subtractive.
An example of a subtractive variety, which is not a Mal'tsev
variety, is the variety of implication algebras \cite{Abb}.

A \textit{subtraction algebra} is a triple $A=(A,-,0)$, where $A$ is
a set, ``$-$\!'' is a binary operation on $A$, and $0$ is a nullary
operation on $A$, satisfying the axioms $x-0=x$ and $x-x=0$. A
\textit{group} can be defined as a subtraction algebra with the
additional axiom $(x-y)-(z-y)=x-z$, and an \textit{abelian group}
can be defined as a subtraction algebra with the stronger axiom
\begin{equation}\label{commutativity of subtraction}(x-y)-(z-t)=(x-z)-(y-t);\end{equation} in both cases,
the group addition $+$ is defined by the equality $x+y=x-(0-y)$,
while the unary operation of taking the inverse $-x$ of an element
$x$, is defined by the equality $-x=0-x$. Note that, to require the
axiom (\ref{commutativity of subtraction}), is the same as to
require that the binary operation $$-:A\x A\ra A$$ is a homomorphism
of subtraction algebras.  Thus, an abelian group is a subtraction
algebra whose algebraic structure is at the same time an internal
subtraction structure in the category $\mathcal{S}$ of all
subtraction algebras. Further, any internal subtraction structure on
an object $A$ in $\mathcal{S}$, coincides with the underlying
subtraction structure of $A$ (see \cite{Bou96},
\cite{BorBou04}\footnote{In \cite{BorBou04}, the term
``protosubtraction'' is used for the operation of subtraction of a
subtraction algebra $A$, and ``subtraction'' is used only when $A$
is a group; in \cite{Bou96} (and also in most cases in
\cite{BorBou04} too), ``$x-y$'' is written as ``$y\setminus x$'',
and accordingly, a ``(proto)subtraction'' is called a
``(proto)devision''.}), and so an internal subtraction algebra in
$\mathcal{S}$ is always an internal abelian group.

The variety $\mathcal{S}$ is clearly a subtractive variety. Further,
the algebraic theory of $\mathcal{S}$ contains a unique nullary
term; equivalently, $\mathcal{S}$ is a \textit{pointed category},
i.e.\! the terminal object in $\mathcal{S}$ coincides with the
initial object. The notion of a subtractive category, introduced in
\cite{Jan05} (see also \cite{Jan06I}), extends the notion of a
pointed subtractive variety to abstract pointed categories --- a
pointed variety is subtractive in the sense of Ursini, if and only
if it is a subtractive category. The purpose of this paper is to
show that for a finitely complete pointed category, with binary
sums, and such that each morphism $f$ in $\C$ can be presented as a
composite $f=me$, where $m$ is a monomorphism and $e$ is an
epimorphism, we have the following result: $\C$ is additive if and
only if both $\C$ and its dual $\C^\mathrm{op}$ are subtractive
(Theorem \ref{thm: The main theorem}). This implies that a Barr
exact \cite{Bar} pointed category $\C$, having binary sums, is an
abelian category, if and only if both $\C$ and $\C^\mathrm{op}$ are
subtractive. In particular, for varieties of universal algebras
(which are always Barr exact and have binary sums), we obtain: a
pointed subtractive variety $\C$ is abelian (i.e.\! $\C$ is the
variety of modules over a fixed ring) if and only if
$\C^\mathrm{op}$ is also subtractive.

\section{The definition of a subtractive category}\label{sec: Preliminaries}

Let $\C$ be a pointed category with binary products. $\C$ is said to
be subtractive, if for any object $A$ in $\C$, any subobject $r:R\ra
A\x A$ of $A\times A$ satisfies the following condition: if there
exist morphisms $f_1, f_2:A\ra R$, such that both diagrams
\begin{equation}\label{dia: Subtractivity A}
\vcenter{ \xymatrix@=3pc{& R\ar[d]^-{r} & & & R\ar[d]^-{r}
\\A\ar[r]_-{(1_A,1_A)}\ar[ur]^-{f_1}& A\x A & & A\ar[r]_-{(1_A,0)}\ar[ur]^-{f_2}& A\x
A}}
\end{equation}
commute, then there exists a morphism $f_3:A\ra R$, such that the
diagram
\begin{equation}\label{dia: Subtractivity B}
\vcenter{\xymatrix@=3pc{& R\ar[d]^-{r}
\\A\ar[r]_-{(0,1_A)}\ar[ur]^-{f_3}& A\x A}}
\end{equation}
commutes. The above condition on the subobject $R$ can be
reformulated, using generalized elements, as follows: if $(a,a)\in
R$ and $(a,0)\in R$ for all $a\in A$, then $(0,a)\in R$ for all
$a\in A$. As we know from \cite{Jan06I}, if $\C$ has finite limits,
then $\C$ is subtractive if and only if the following condition is
satisfied for any object $A$ in $\C$, and any subobject $R$ of
$A\times A$: for all $a\in A$,\begin{equation}
\label{implication I}(a,a)\in R\;\wedge\; (a,0)\in R\;\;\Rightarrow\;\; (0,a)\in
R.\end{equation} The condition (\ref{implication I}) on $R$ states
precisely that $R$ is (regarded as a binary relation on $A$)
\textit{$M$-closed} in the sense of
\cite{Jan06I}, where $M$ is the extended matrix
$$M=\left(\begin{array} {cc|c} x & x & 0 \\ x & 0 & x
\end{array}\right)$$
of terms in the algebraic theory $\T$ of pointed sets --- $0$ is the
unique constant of $\T$, while $x$ is an arbitrary variable of $\T$.
Note that this matrix directly arises from the system of equations
$$
\left\{\begin{array}{l} x - x = 0,\\ x - 0 = 0,\end{array}\right.
$$
which defines an operation of subtraction. In the present paper we
write such matrices mostly in the transposed form $$M=\left(\begin{array}{cc}x & x\\
x & 0 \\ \hline 0 & x\end{array}\right).$$ \textit{$R$ is
$M$-closed} means that for any
\textit{interpretation} $$M^a=\left(\begin{array}{cc}a & a\\
a & 0 \\ \hline 0 & a\end{array}\right)$$ of $M$, where $a\in A$ and
$0$ is the base point of $A$, if the first two rows of $M^a$ belong
to $R$ (considered as generalized elements of $A\x A$), then also
the third row (i.e.\! the row below the horizontal line) belongs to
$R$. Thus, $R$ is $M$-closed if and only if the implication
(\ref{implication I}) is satisfied for all $a\in A$, as already said
above. Recall from \cite{Jan06I}, that if we take $$M=\left(\begin{array}{cc}x & y\\
x & 0 \\ \hline 0 & y\end{array}\right),$$ then, again, a (finitely
complete pointed) category with $M$-closed relations is the same as
a (finitely complete) subtractive category. If we further
modify $M$ by adding to it columns of the form $$\left(\begin{array}{c}u\\
u\\ \hline 0\end{array}\right),\textrm{ or of the form }\left(\begin{array}{c}v\\
0\\ \hline v\end{array}\right),$$ then again (by Proposition 1.7 of
\cite{Jan06II}), a subtractive category is the same as a pointed
category in which every $n$-ary relation $R\ra A^n$ is $M$-closed,
where now $M$ is the new matrix, and $n$ denotes the number of
columns of $M$.

\section{\textit{Additive} is the same as \textit{subtractive and cosubtractive}}
\label{sec: Additivity via subtractivity}

A \textit{semi-abelian category} in the sense of G.\,Janelidze,
L.\,M\'arki, and W.\,Tholen \cite{JanMarTho02}, is a pointed
category $\C$ satisfying the following conditions:
\begin{itemize}
\item $\C$ is Barr exact \cite{Bar},

\item $\C$ has binary sums,

\item $\C$ is protomodular in the sense of D.\,Bourn
\cite{Bou91}.
\end{itemize}
As shown in \cite{BouJan03}, a pointed variety of universal algebras
is semi-abelian if and only if it is \textit{classically ideal
determined} in the sense of A.\,Ursini \cite{Urs94} (this is the
same as \textit{BIT speciale} in the sense of \cite{Urs73}).

As observed in \cite{Jan05}, any semi-abelian category is
subtractive.

The following equation for categories was obtained in
\cite{JanMarTho02}:
$$\textrm{semi-abelian }+\textrm{ semi-abelian}^\mathrm{op}=\textrm{abelian}.$$
This equation says: a category $\C$ is an abelian category if and
only if both $\C$ and its dual category $\C^\mathrm{op}$ are
semi-abelian. From Theorem \ref{thm: The main theorem} below, it
follows that, more generally, for any Barr exact category with
binary sums, we have (see Corollary \ref{cor: Abelian categories},
and see also Remark \ref{rem: Strongly unital categories}):
$$\textrm{subtractive }+\textrm{ subtractive}^\mathrm{op}=\textrm{abelian}.$$
Recall that a Barr exact category $\C$ is an abelian category if and
only if it is additive, i.e.\! if and only if $\C$ is enriched in
the category of abelian groups --- each object in $\C$ is equipped
with an internal abelian group structure, so that each morphism in
$\C$ is a homomorphism of internal abelian groups. More generally,
we have:

\begin{lemma}\label{lem: additivity via subtraction algebras}
A category with finite products is additive if and only if each
object in it is equipped with an internal subtraction structure, so
that each morphism is a homomorphism of internal subtraction
algebras.
\end{lemma}

\begin{proof}
This follows immediately from the definition of an additive
category, and the known fact that if the operation of subtraction of
a subtraction algebra is homomorphic, then the subtraction algebra
becomes an abelian group.
\end{proof}

\begin{theorem}\label{thm: The main theorem}
Let $\C$ be a finitely complete pointed category with binary sums,
in which

\begin{itemize}
\item[(*)] any morphism $f:X\ra Y$ can be decomposed $f=me$ into a
monomorphism $m$ and an epimorphism $e$.
\end{itemize}

\noindent Then, $\C$ is an additive category if and only if both
$\C$ and $\C^\mathrm{op}$ are subtractive.
\end{theorem}

\begin{proof}
If $\C$ is additive, then it is subtractive. Indeed, then for any
pair of commutative diagrams (\ref{dia: Subtractivity A}), we can
form the third one (\ref{dia: Subtractivity B}), by taking in it
$f_3=f_1-f_2$, where ``$-$'' denotes the operation of subtraction on
$\mathrm{hom}(A,R)$, induced by the operation of subtraction of the
abelian group structure of $R$. We will then have
$$rf_3=r(f_1-f_2)=rf_1-rf_2=(1_A,1_A)-(1_A,0)=(1_A-1_A,1_A-0)=(0,1_A),$$
so the diagram (\ref{dia: Subtractivity B}) indeed commutes. On the
other hand, if $\C$ is additive, then so is $\C^\mathrm{op}$, and
hence $\C^\mathrm{op}$ is subtractive. This proves the first part of
the theorem.

Assume now that both $\C$ and $\C^\mathrm{op}$ are subtractive.
According to Lemma \ref{lem: additivity via subtraction algebras},
to show that $\C$ is additive, it suffices to show that each object
$A$ in $\C$ is equipped with an internal subtraction structure,
$A=(A,-,0)$, so that each morphism $f:A\ra A'$ in $\C$ is a
homomorphism of (internal) subtraction algebras. Let $A$ be an
object in $\C$. Then we can form a commutative diagram
$$\xymatrix@=3pc{
 & R\ar[dr]^-{r=(r_1,r_2)} &\\
A+A \ar[ur]^-{q=\columnIIs{q_1}{q_2}}\ar[rr]_-{\matrixIIxIIs
{1_A}{1_A}{1_A}{0}} & & A\x A\\}
$$
where $r$ is a monomorphism and $q$ is an epimorphism. This
commutative diagram yields the following two commutative diagrams:
$$
\xymatrix@=3pc{
  R\ar[dr]^-{r_1}&\\
A+A \ar[u]^-{q}\ar[r]_-{\columnIIs{1_A}{1_A}} & A\\}
\spaces\spaces\spaces
\xymatrix@=3pc{ R\ar[dr]^-{r_2} &\\
A+A \ar[u]^-{q}\ar[r]_-{\columnIIs{1_A}{0}} & A\\}
$$
Since $\C^\mathrm{op}$ is subtractive, there exists a morphism
$s:R\ra A$ such that the diagram
$$
\xymatrix@=3pc{
R\ar[dr]^-{s} &\\
A+A \ar[u]^-{q}\ar[r]_-{\columnIIs{0}{1_A}} & A\\}
$$
commutes. Next, we regard $R$ as the subobject $r:R\ra A\x A$ of
$A\x A$. Then $s:R\ra A$ becomes a partial binary operation on $A$.
Commutativity of the two diagrams
$$
\xymatrix@=3pc{
 & R\ar[d]^-{r}\\
A \ar[ur]^-{q_1}\ar[r]_-{(1_A,1_A)} & A\x A\\} \spaces\spaces\spaces
\xymatrix@=3pc{ & R\ar[d]^-{r}\\
A \ar[ur]^-{q_2}\ar[r]_-{(1_A,0)} & A\x A\\}
$$
and the fact that $\C$ is subtractive, yield that there exists a
morphism $p:A\ra R$ such that the diagram
$$
\xymatrix@=3pc{
 & R\ar[d]^-{r}\\
A \ar[ur]^-{p}\ar[r]_-{(0,1_A)} & A\x A\\}
$$
commutes. From the commutativity of the last four diagrams we
obtain: For all $a\in A$,
$$(a,a),(a,0),(0,a)\in R$$ and, further,
$s(a,a)=0$ and $s(a,0)=a$. We now show $R=A\x A$. Consider the
ternary relation $R'\ra A^3$ defined as follows: $$(a,b,c)\in
R'\;\;\Leftrightarrow\;\; (b,c)\in R\;\wedge\;(a,s(b,c))\in R.
$$
Since $\C$ is subtractive, $R'$ is closed with respect to the matrix
$$\left(\begin{array}{ccc}x & y & y \\ 0 & 0 &
y\\\hline x & y & 0\\\end{array}\right).$$ For each $a,b\in A$ we
then have:
$$\begin{array}{rlcc}
\smallskip (b,b)\in R\;\wedge\;(a,s(b,b))=(a,0)\in R & \Rightarrow & (a,b,b)\in R', & \\
(0,b)\in R\;\wedge\;(0,s(0,b))\in R & \Rightarrow & (0,0,b)\in R', & \\
\smallskip & &\Downarrow &\\
(a,b)=(a,s(b,0))\in R & \Leftarrow & (a,b,0)\in R'. &\\
\end{array}$$
Thus, for all $a,b\in A$, $(a,b)\in R$, which shows $R=A\x A$. So
$s$ is a full binary operation on $A$, and since we already know
that it satisfies the identities $s(a,a)=0$ and $s(a,0)=a$, we get
that $(A,s,0)$ is a subtraction algebra. We have thus shown
\begin{itemize}
\item that the morphism $r$ is an isomorphism; this implies that the
morphism $$\xymatrix@=3pc{ A+A\ar[rr]_-{\matrixIIxIIs
{1_A}{1_A}{1_A}{0}} & & A\x A\\}
$$ is an epimorphism;

\item that there exists a subtraction algebra $(A,-,0)$; note that the
subtraction axioms for the operation $-:A\x A\ra A$ state precisely
that the diagram
$$\xymatrix@=3pc{
A+A \ar@/_-2pc/[rrr]^-{\columnIIs{0}{1_A}}\ar[rr]_-{\matrixIIxIIs
{1_A}{1_A}{1_A}{0}} & & A\x A\ar[r]_-{-} & A\\}
$$
commutes.
\end{itemize}
Now let $(A',-',0')$ be another subtraction algebra, obtained from
an object $A'$, like $(A,-,0)$ is obtained from $A$. Let $f$ be any morphism
$f:A\ra A'$. In the diagram
$$
\xymatrix@=3pc{ A+A\ar[d]_-{f+f}\ar[rr]^-{\matrixIIxIIs
{1_A}{1_A}{1_A}{0}} & &
A\x A\ar[d]_-{f\x f}\ar[r]^-{-} & A\ar[d]^-{f}\\
A'+A'\ar[rr]_-{\matrixIIxIIs {1_{A'}}{1_{A'}}{1_{A'}}{0'}} & & A'\x
A'\ar[r]_-{-'} & A'\\}
$$
the left inner rectangle commutes and the outer rectangle commutes.
Since the first morphism in the top row is an epimorphism, this
implies that the right inner rectangle commutes. We thus obtain that
$f$ is a homomorphism $f:(A,-,0)\ra (A',-',0')$ of subtraction
algebras. This concludes the proof.
\end{proof}

\begin{corollary}\label{cor: Abelian categories}
A Barr exact category $\C$ with binary sums is an abelian category,
if and only if both $\C$ and its dual $\C^\mathrm{op}$ are
subtractive.
\end{corollary}

\begin{proof}
An ablian category is the same as an additive Barr exact category.
On the other hand, a Barr exact category, being a regular category,
has finite limits and any morphism $f$ in it can be decomposed
$f=me$ into a monomorphism $m$ and an epimorphism $e$ (moreover, any
morphism can be decomposed into a monomorphism and a regular
epimorphism). After these observations, it is only left to apply
Theorem \ref{thm: The main theorem}.
\end{proof}

\begin{remark}\label{rem: Strongly unital categories}
I do not know if Theorem \ref{thm: The main theorem} remains true if
$\C$ in it does not have the property (*). However, it does remain
true if in addition we replace ``subtractive'' in the theorem with
``strongly unital'' in the sense of D.~Bourn \cite{Bou96}; this is
an immediate consequence of the following known facts (see
\cite{BorBou04} and \cite{Jan05}):
\begin{itemize}
\item[(i)] a finitely complete category is additive if and only if it
is subtractive and enriched in the category of commutative monoids,

\item[(ii)] a category $\C$ with binary sums and products is enriched in
the category of commutative monoids if and only if both $\C$ and its
dual $\C^\mathrm{op}$ are unital in the sense of D.~Bourn
\cite{Bou96},

\item[(iii)] a finitely complete category is strongly unital if and only if
it is subtractive and unital.
\end{itemize}
\end{remark}


\begin{thebibliography}{...}
\bibitem{Abb}
J. C. Abbot,
Implication algebras,
\textit{Bull. Math. Soc. Sci. Roumanie} 11(59)
1967, 3-23 (1968).

\bibitem{Bar}
M. Barr, \textit{Exact Categories}, Lecture Notes in Mathematics
236, Springer, Berlin, 1971, 1-120.

\bibitem{BorBou04}
F. Borceux and D. Bourn,
\textit{Mal'cev, protomodular, homological and semi-abelian categories},
Mathematics and its Applications 566, Kluwer,
2004.

\bibitem{Bou91}
D. Bourn,
Normalization equivalence, kernel equivalence and affine categories,
\textit{Springer Lecture Notes in Mathematics} 1488,
1991, 43-62.

\bibitem{Bou96}
D. Bourn,
Mal'cev categories and fibration of pointed objects,
\textit{Applied Categorical structures} 4,
1996, 307-327.

\bibitem{BouJan03}
D. Bourn and G. Janelidze, Characterization of protomodular
varieties of universal algebras, \textit{Theory and Applications of
Categories}, Vol. 11, No. 6, 2003, 143-447.

\bibitem{JanMarTho02}
G. Janelidze, L. M\'arki, and W. Tholen,
Semi-abelian categories,
\textit{Journal of Pure and Applied Algebra} 168,
2002, 367-386.

\bibitem{Jan05}
Z. Janelidze,
Subtractive categories,
\textit{Applied Categorical Structures}, Vol. 13, No. 4,
2005, 343-350.

\bibitem{Jan06I}
Z. Janelidze,
Closedness properties of internal relations I: A unified approach to Mal'tsev, unital and subtractive categories,
\textit{Theory and Applications of Categories}, Vol. 16, No. 12,
2006, 236-261.

\bibitem{Jan06II}
Z. Janelidze,
Closedness properties of internal relations II: Bourn localization,
\textit{Theory and Applications of Categories}, Vol. 16, No. 13,
2006, 262-282.

\bibitem{Jan06III}
Z. Janelidze, Closedness properties of internal relations III:
Pointed protomodular categories, submitted for publication.

\bibitem{Mal54}
A. I. Mal'tsev, On the general theory of algebraic systems,
\textit{USSR Mathematical Sbornik N.S.} 35, 1954, 3-20.

\bibitem{Urs73}
A. Ursini,
Osservazioni sulla variet\'a BIT,
\textit{Bolletino della Unione Matematica Italiana} 8,
1973, 205-211.

\bibitem{Urs94}
A. Ursini,
On subtractive varieties, I,
\textit{Algebra Universalis} 31,
1994, 204-222.
\end{thebibliography}
\end{document}